\documentclass[12pt,twoside]{article}

\usepackage{amsmath}
\usepackage{amsfonts}
\usepackage{amsthm,amssymb}
\usepackage{color}
\newcommand{\RR}{\mathbb R}

\newcommand{\email}[1]{{\small E-mail: {\textsf {#1}}}}
\newtheorem{defi}{Definition}%[section]
\newtheorem{pr}[defi]{Proposition}
\newtheorem{theo}[defi]{Theorem}
\newtheorem{Le}[defi]{Lemma}

\newtheorem{Rem}[defi]{Remark}

\begin{document}
\title{Looking for critical nonlinearity in the one-dimensional quasilinear Smoluchowski-Poisson system} 

\author{Tomasz Cie\'slak\footnote{Institute of Applied Mathematics, Warsaw University, Banacha 2, 02-097 Warszawa, Poland. \email{T.Cieslak@impan.gov.pl}} \kern8pt \& \kern8pt
Philippe Lauren\c cot\footnote{Institut de Math\'ematiques de Toulouse, CNRS UMR~5219, Universit\'e de Toulouse, F--31062 Toulouse cedex 9, France. \email{Philippe.Laurencot@math.univ-toulouse.fr}}}
\date{\today}
\maketitle

\begin{abstract}
It is known that classical solutions to the one-dimensional quasilinear Smoluchowski-Poisson system with nonlinear diffusion $a(u)=(1+u)^{-p}$ may blow up in finite time if $p>1$ and exist globally if $p<1$. The case $p=1$ thus appears to be critical but it turns out that all solutions are global also in that case. Two classes of diffusion coefficients are actually identified in this paper, one for which all solutions to the corresponding quasilinear Smoluchowski-Poisson system are global and the other one leading to finite time blow-up for sufficiently concentrated initial data. The cornerstone of the proof are an alternative formulation of the Smoluchowski-Poisson system which relies on a novel change of variables and a virial identity.    
\end{abstract}

%%%%%%%%%%%%%%%%%%%%%%%%%%%%%%%%%%%%%%%%%%%%%%%
%%%%%%%%%%%%%%%%%%%%%%%%%%%%%%%%%%%%%%%%%%%%%%%
\section{Introduction.}\label{intro}
%%%%%%%%%%%%%%%%%%%%%%%%%%%%%%%%%%%%%%%%%%%%%%%
%%%%%%%%%%%%%%%%%%%%%%%%%%%%%%%%%%%%%%%%%%%%%%%

In one space dimension, the quasilinear Smoluchowski-Poisson (SP) system reads
\begin{eqnarray}
\label{he1}
\partial_t u &=& \partial_x \left( a(u)\ \partial_x u - u\ \partial_x v \right) \;\;\mbox{in}\;\;(0,T)\times (0,1),\\
\label{he12}
0&=& \partial_{x}^2 v + u - M \;\;\mbox{in}\;\;(0,T)\times (0,1),\\
\label{he2}
\partial_x A(u)(t,0) &=& \partial_x A(u)(t,1)= \partial_x v(t,0) = \partial_x v(t,1) = 0\;\;\mbox{ for } \;\; t\in (0,T)\\
\label{he3}
u(0)&=& u_0\geq 0 \;\;\mbox{in} \;\;(0,1),\;\; \langle v(t) \rangle = 0\;\;\mbox{for}\;\; t\in (0,T),
\end{eqnarray}
where $(u,v):(0,T)\times (0,1)\rightarrow \RR$ are the unknown functions, $a\in C^1((0,\infty))$ is a positive function, $A$ is an indefinite integral of $a$, and the constant $M:=\langle u_0\rangle$, where $\langle f\rangle := \int_I f(x)\ dx /|I|$ denotes the mean value of a given function $f\in L^1(I)$ defined on an interval $I$ of $\RR$ of length $|I|$.  
Such systems appear in the modelling of self-gravitating systems in astrophysics \cite{Ch03} and chemotactic processes in biology \cite{KS70}. One particular feature of this system is that there is a competition between the diffusion $\partial_x(a(u)\partial_x u)$ and the drift term $\partial_x(u\partial_x v)$ that may result in the occurrence of blow-up in finite time, a phenomenon related to the gravitational/chemotactic collapse in astrophysics and biology. Roughly speaking, it has been shown in any space dimension $d\ge 1$ that, if $a(r)\le C\ (1+r)^{-p}$ for $r>0$ and some $p>(2-d)/d$, there are solutions to the quasilinear SP system blowing up in finite time while all solutions are global if $a(r)\ge C\ (1+r)^{-p}$ for $r>0$ and some $p<(2-d)/d$ \cite{CW08} (see also \cite{ChS04,Su07} for a related system in $\RR^d$). Thus, the exponent $p=(2-d)/d$ appears to be a critical exponent as it corresponds to the borderline growth for the nonlinear diffusion coefficient that separates two different behaviours. Furthermore, if $d\ge 2$, it has been shown in \cite{cl1,GZ98,Na95} that, if $a(r) = (1+r)^{(d-2)/d}$, there is a threshold value of $M$ above which finite time blow-up takes place for some initial data while global existence is guaranteed if $M$ is sufficiently small (see also \cite{BCLxx,ChS04,DP04,Su07} and the references therein for the quasilinear SP system in $\RR^d$). The one-dimensional case $d=1$ (corresponding to $a(r)=1/(1+r)$) is not covered by the above mentioned results and the main purpose of this paper is to elucidate what happens in that particular case.

As we shall see below, it turns out that the situation for $a(r)=1/(1+r)$ in one space dimension is strikingly different from that  encountered in higher space dimensions and that all solutions are global. We actually identify a class of diffusion coefficients   encompassing the one found in \cite{CW08} for which all solutions to \eqref{he1}-\eqref{he3} are global. Our approach relies on a new formulation of the SP system \eqref{he1}-\eqref{he3} which transforms it into a single quasilinear parabolic equation. While it is well-known that the indefinite integral (with respect to space) of $u$ vanishing at zero solves a quasilinear parabolic equation \cite{CW08,JL92}, we go one step beyond and derive the equation satisfied by the first derivative of the inverse of  the indefinite integral of $u$. The latter is the solution of a nonlinear reaction-diffusion equation with homogeneous Neumann boundary conditions, the reaction term being affine and the diffusion term nonlinear. This transformation is described in details in Section~\ref{cv}. Let us emphasize at this point that the blow-up in finite time of a solution to \eqref{he1}-\eqref{he3} amounts to the so-called ``touch-down'' of the solution to the transformed equation: in other words, a solution to \eqref{he1}-\eqref{he3} blows up in finite time if and only if the corresponding solution to the transformed equation reaches zero in finite time.  Section~\ref{wp} is then devoted to the local-in-time existence and uniqueness of solutions. At this point we take advantage of the transformation introduced in Section~\ref{cv} to provide an simple proof of the uniqueness which is somehow related to the recent use of Wasserstein distances to study nonlinear diffusion equations \cite{CF05,CGT04,DFW08}. In Section~\ref{glex}, we show that, for a suitable class of functions $a$,  ``touch-down'' cannot occur in finite time for the transformed problem and as a corollary we obtain existence of global  solutions to \eqref{he1}-\eqref{he3}. Furthermore, as a byproduct of our analysis, we also establish that the solutions are bounded in $L^\infty$ uniformly with respect to time. This result applies in particular to
\begin{equation}
\label{exglex}
a(r)=\frac{1}{(1+r)^p}\,, \quad p\in (-\infty,1]\,, \;\;\mbox{ and }\;\; a(r) = \frac{1}{(1+r) (\log{(1+r)})^\alpha}\,, \quad \alpha\in (-\infty,1]\,.
\end{equation}
In the subsequent Section~\ref{ftbu} we show that, whatever the value of the mass $M$ is, singularities of solutions appear in finite time for some diffusion coefficients between $(1+r)^{-1}$ and  $(1+r)^{-p}, p>1$, namely
\begin{equation}
\label{exftbu} 
a(r) \le \frac{1}{(1+r)^p}\,, \quad p\in (1,2]\,, \;\;\mbox{ and }\;\; a(r) \le \frac{1}{(1+r)\ \left(\log(1+r)\right)^{\alpha}}\,, \quad \alpha\in (1,2]\,. 
\end{equation}
Thus we eliminate some candidates for the critical nonlinearity (if any!) for \eqref{he1}-\eqref{he3}. Observe that the gap between the functions listed in \eqref{exglex} and in \eqref{exftbu} is quite narrow. The proof of the finite time blow-up relies on some virial identities which differ from the ones used originally for the classical SP system corresponding to $a(r)=1$ \cite{BN94,Na95}. Here again, while the virial identity is written in Section~\ref{ftbu} in terms on the indefinite integral of $u$, we mention that the initial guess for this identity was made on the transformed equation introduced in Section~\ref{cv}, thus emphasizing once more the usefulness of this alternative formulation. Virial identities of the same type have recently allowed us to obtain some finite time blow-up results for the parabolic-parabolic quasilinear one-dimensional Keller-Segel system \cite{cl2} (in that case, the equation \eqref{he12} for $v$ is replaced by a parabolic equation). Up to our knowledge this is the only result concerning blow-up in finite time for the parabolic-parabolic Keller-Segel system besides the remarkable (and more precise) result in \cite{HV97}. 

%%%%%%%%%%%%%%%%%%%%%%%%%%%%%%%%%%%%%%%%%%%%%%%
%%%%%%%%%%%%%%%%%%%%%%%%%%%%%%%%%%%%%%%%%%%%%%%
\section{Change of variables. Blow-up means touch-down.}\label{cv} 
%%%%%%%%%%%%%%%%%%%%%%%%%%%%%%%%%%%%%%%%%%%%%%%
%%%%%%%%%%%%%%%%%%%%%%%%%%%%%%%%%%%%%%%%%%%%%%%

In this section we introduce the change of variables that enables us to reformulate the quasilinear SP system \eqref{he1}-\eqref{he3} in such a way that we arrive at a one-dimensional quasilinear parabolic equation with an affine right-hand side. It involves the space derivative $f$ of the inverse of the cumulative distribution of $u$. Its surprising  applicability to \eqref{he1}-\eqref{he3} has its origins in the fact that it transforms the drift term $\partial_x\left( u\ \partial_x v \right)$ in a much simpler one, see \eqref{main} below. To the best of our knowledge, such a change of variables has not appeared so far in the literature in the study of global existence of solutions. However, since the $q$-Wasserstein distance between two functions is related to the $L^q$-norm of the difference of the inverses of their indefinite integrals in one space dimension (see, e.g.,  \cite{Vi03}), our approach is somehow related to some recent works on degenerate diffusion equations where Wasserstein metrics are used to investigate the large time behaviour and the speed of propagation of solutions to these equations (see, e.g., \cite{CF05,CGT04,DFW08,Vi03} and the references therein).

\medskip

Let $T>0$ and assume that $u\in C([0,T)\times [0,1])\cap C^{1,2}((0,T)\times[0,1])$ is a classical solution to \eqref{he1}-\eqref{he3} such that  $u>0$ in $[0,T)\times(0,1)$. For $t\in [0,T)$, we define the indefinite integral of $u(t)$ vanishing at zero by 
\begin{equation}\label{indefint}
U(t,x) :=\int_0^x u(t, z)dz\,, \qquad x\in [0,1]\,,
\end{equation}
and observe that the positivity of $u$ warrants that $x\mapsto U(t,x)$ is an increasing function from $[0,1]$ onto $[0,M]$. For each $t\in [0,T)$, let $y\mapsto F(t,y)$ denote its inverse function which is well-defined, maps $[0,M]$ onto $[0,1]$, and satisfies
\begin{equation}\label{zero}
U(t,F(t,y))=y\,, \qquad (t,y)\in [0,T)\times [0,M]\,. 
\end{equation}
We next set $f(t,y):=\partial_y F(t,y)$ for $(t,y)\in [0,T)\times (0,M)$. Differentiating \eqref{zero} once with respect to $y$, we first see  that 
\begin{equation}\label{odw}
f(t,y)\ u(t,F(t,y))=1\,, \qquad (t,y)\in [0,T)\times [0,M]\,.
\end{equation} 
Next, differentiating \eqref{odw}  with respect to the space variable and differentiating \eqref{zero} with respect to the time variable we arrive at 
\begin{equation}\label{dwa}
\partial_x u(t,F(t,y))\ f(t,y)^2+u(t,F(t,y))\ \partial_y f(t,y)=0\,, \qquad (t,y)\in (0,T)\times (0,M)\,,
\end{equation}
and
\begin{equation}\label{raz}
\partial_t U(t,F(t,y))+u(t,F(t,y))\ \partial_tF(t,y)=0\,, \qquad (t,y)\in (0,T)\times (0,M)\,.
\end{equation}
Now, integrating \eqref{he1} and \eqref{he2} with respect to space over $(0,x)$ we obtain
\[
\partial_t U(t,x)=a(u(t,F(t,y))\ \partial_x u(t,x) - u(t,x)\ \partial_x v(t,x),\;\;\;-\partial_x v(t,x)=U(t,x)-Mx,
\]
for $(t,x)\in (0,T)\times (0,1)$, so that 
\begin{equation}\label{B}
\partial_t U(t,x)=a(u(t,x))\ \partial_x u(t,x) + u(t,x)\ (U(t,x)-Mx)\,, \qquad (t,x)\in (0,T)\times (0,1)\,.
\end{equation}
Therefore, for $(t,y)\in (0,T)\times (0,M)$, we have 
\begin{equation}\label{trzy}
\partial_t U(t,F(t,y))=a(u(t,F(t,y))\ \partial_x u(t,F(t,y)) + u(t,F(t,y))\ (y-MF(t,y))
\end{equation}
by \eqref{zero}. We next use \eqref{odw}, \eqref{dwa}, and \eqref{raz} along with the positivity of $u$ to deduce that
\[
-u(t,F(t,y))\ \partial_tF(t,y)=-\frac{u(t,F(t,y))a(u(t,F(t,y)))\partial_yf(t,y)}{f(t,y)^2}+u(t,F(t,y))(y-MF(t,y))
\]
and thus
\[
\partial_t F(t,y)=\frac{1}{f(t,y)^2}\ a\left( \frac{1}{f(t,y)} \right)\ \partial_yf(t,y)-y+MF(t,y)\,, \qquad (t,y)\in (0,T)\times (0,M)\,.
\] 
Differentiating both sides of the previous identity with respect to $y$ and introducing
\begin{equation}\label{Psi}
\Psi'(r):=\frac{1}{r^2}\ a\left( \frac{1}{r} \right) \;\;\mbox{for any}\;\;r>0\,, \qquad \Psi(1):=0\,,
\end{equation} 
we realize that $f$ solves
\begin{equation}\label{main}
\partial_t f=\partial_y^2 \Psi(f)-1+Mf\,, \qquad (t,y)\in (0,T)\times (0,M)\,.
\end{equation} 
Next, recalling \eqref{dwa}, the boundary conditions \eqref{he3} become
\begin{equation}\label{boundary}
\partial_y f(t,0)=\partial_y f(t,M)=0\,, \qquad t\in (0,T)\,,
\end{equation}
while the initial value for $f$ is given by \eqref{he3}, \eqref{zero}, and \eqref{odw}, and reads 
\begin{equation}\label{init}
f(0,y)=f_0(y):= \frac{1}{u_0(F(0,y))}\,, \qquad y\in (0,M)\,.
\end{equation}
Moreover the conservation of mass yields
\begin{equation}\label{mas}
\int_0^M f(t,y) dy=F(t,M)-F(t,0)=1\,, \qquad t\in (0,T)\,.
\end{equation}

We have thus established the following result:
\begin{pr}\label{altform}
Let $T>0$ and consider a classical solution $u\in C([0,T)\times [0,1])\cap C^{1,2}((0,T)\times[0,1])$ to \eqref{he1}-\eqref{he3} such that  $u>0$ in $[0,T)\times(0,1)$. Defining $U$and $F$ by \eqref{indefint} and \eqref{zero}, respectively, the function $f\in C([0,T)\times [0,M])\cap C^{1,2}((0,T)\times[0,M])$ defined by $f:=\partial_y F$ is a classical solution to \eqref{main}-\eqref{init} satisfying the conservation  of mass \eqref{mas}.
\end{pr}

\begin{Rem}\label{td}
Since the main issue of our analysis is whether solutions to \eqref{he1}-\eqref{he3} exist globally or blow up in finite time, we emphasize here that, according to \eqref{odw}, finite time blow-up at $T_{max}$ of a solution $u$ to \eqref{he1}-\eqref{he3} corresponds to finite time ``touch-down'' at $T_{max}$ for the related solution $f$ to \eqref{main}-\eqref{init}, namely,
\begin{equation}\label{tdf}
\lim_{t\to T_{max}}\ \min_{y\in [0,M]}{ \{f(t,y)\} } = 0\,. 
\end{equation}
In other words, $f$ ceases to be positive at time $T_{max}$ and reaches the zero value somewhere in the interval $[0,M]$.
\end{Rem}

%%%%%%%%%%%%%%%%%%%%%%%%%%%%%%%%%%%%%%%%%%%%%%%
%%%%%%%%%%%%%%%%%%%%%%%%%%%%%%%%%%%%%%%%%%%%%%%
\section{Well-posedness.}\label{wp}
%%%%%%%%%%%%%%%%%%%%%%%%%%%%%%%%%%%%%%%%%%%%%%%
%%%%%%%%%%%%%%%%%%%%%%%%%%%%%%%%%%%%%%%%%%%%%%%

We first recall the following result \cite[Theorem 1.3]{CW08}. 
\begin{pr}\label{local1}
Consider a positive function $a\in C^1([0,\infty))$ and a nonnegative initial condition $u_0\in C([0,1])$ such that $\langle u_0\rangle=M$ for some $M>0$. Then there exist a time $T_{max}\in (0,\infty]$ and a unique solution 
$$
(u,v)\in C([0,T_{max})\times[0,1];\RR^2)\cap C^{1,2}((0,T_{max})\times [0,1];\RR^2)\,, \qquad u\ge 0\,,
$$ 
to \eqref{he1}-\eqref{he3} such that $\langle u(t)\rangle=M$ for any $t\in [0,T_{max})$. Moreover, if $T_{max}<\infty$
then $\left\|u(t)\right\|_\infty \longrightarrow \infty$ as $t\rightarrow T_{max}$. By the strong maximum
principle we also have $u>0$ in $(0,T_{max})\times (0,1)$.
\end{pr}

Proposition~\ref{local1} does not apply to the case $a(r)=1/r$, $r>0$, because of the singularity of $a$ at zero. In that case, we nevertheless have the following result under the stronger assumption that the minimum of $u_0$ is positive.

\begin{pr}\label{local2}
Consider a positive function $a\in C^1((0,\infty))$ and an initial condition $u_0\in C([0,1])$ such that $u_0\ge m_0>0$ and $\langle u_0\rangle=M$ for some $M>0$ and $m_0\in (0,M)$. Then there exist a time $T_{max}\in (0,\infty]$ and a unique solution
$$
(u,v)\in C([0,T_{max})\times[0,1];\RR^2)\cap C^{1,2}((0,T_{max})\times [0,1];\RR^2)\,,
$$ 
to \eqref{he1}-\eqref{he3} such that 
\begin{equation}\label{swing}
\langle u(t)\rangle=M \;\;\mbox{ and }\;\; u(t,x)\ge \frac{M m_0}{m_0+e^{Mt} (M-m_0)}>0 \;\;\mbox{ for all }\;\; (t,x)\in [0,T_{max})\times [0,1]\,.
\end{equation}
In addition, $T_{max}<\infty$ implies that $\left\|u(t)\right\|_\infty \longrightarrow \infty$ as $t\rightarrow T_{max}$.
\end{pr}

\noindent{\bf Proof.}  
We first prove the existence of a solution by a standard approximation argument: more precisely, for $\varepsilon\in (0,1)$ and $r\ge 0$, we put $a_\varepsilon(r):=a(r+\varepsilon)$. As $a_\varepsilon$ is a positive function in $C^1([0,\infty))$, we infer from Proposition~\ref{local1} that there exists a unique classical solution $(u^\varepsilon, v^\varepsilon)$ to \eqref{he1}-\eqref{he3} with $a_\varepsilon$ instead of $a$ which is defined on its maximal existence time interval $[0,T_{max}^\varepsilon)$. Moreover by standard parabolic regularity results \cite{LSU} we have 
\begin{equation}\label{reg}
\|u^\varepsilon\|_{C^{1+\alpha,2+2\alpha}((0,T)\times[0,1])}\leq C(T)\ \|u^\varepsilon\|_{L^\infty((0,T)\times (0,1))}
\end{equation}
for some $\alpha\in (0,1/2)$ and $T\in (0,T_{max}^\varepsilon)$. By (\ref{he1})-(\ref{he2}) $u^\varepsilon$ satisfies
\begin{eqnarray}
\partial_t u^\varepsilon & = & \partial_x \left( a\left(u^\varepsilon \right)\ \partial_x u^\varepsilon \right) - \partial_x v^\varepsilon\ \partial_x u^\varepsilon + (u^\varepsilon)^2 - u^\varepsilon M\,, \qquad (t,x)\in \left( 0,T_{max}^\varepsilon \right)\times (0,1)\,, \label{aux}\\
& & \partial_x u^\varepsilon(t,0) = \partial_x u^\varepsilon(t,1)=0\,, \qquad t\in \left( 0,T_{max}^\varepsilon \right)\,. \label{auxb}
\end{eqnarray} 
Introducing
$$
\Sigma(t):= \frac{\|u_0\|_\infty}{1-\|u_0\|_\infty\ t} \;\;\mbox{ and }\;\; \sigma(t) := \frac{M m_0}{m_0+e^{Mt} (M-m_0)}\,, \qquad t\in \left( 0, \min{\left\{ T_{max}^\varepsilon , \|u_0\|_\infty^{-1} \right\}} \right)\,,
$$
we notice that $\Sigma$ is a supersolution to \eqref{aux}-\eqref{auxb} while $\sigma$ is a subsolution to \eqref{aux}-\eqref{auxb} with 
$$
\sigma(0)=m_0 \le u_0(x) \le \|u_0\|_\infty = \Sigma(0)\,, \qquad x\in [0,1]\,.
$$
The comparison principle then warrants that
\begin{equation}\label{encadr}
\sigma(t) \le u^\varepsilon(t,x) \le \Sigma(t)\,, \qquad (t,x) \in \left( 0, \min{\left\{ T_{max}^\varepsilon , \|u_0\|_\infty^{-1} \right\}} \right)\times (0,1)\,.
\end{equation}
According to Proposition~\ref{local1}, a first consequence of \eqref{encadr} is that $T_{max}^\varepsilon \ge \|u_0\|_\infty^{-1}$. We next infer from \eqref{reg}, \eqref{encadr}, and the Arzel\`a-Ascoli theorem that there are a nonnegative function $u$ and a subsequence of $(u^\varepsilon)$ (not relabeled) such that 
\begin{equation}
\label{cvg}
u^\varepsilon \longrightarrow u \;\;\mbox{ in }\;\; C\left( \left[ 0,\|u_0\|_\infty^{-1} \right)\times [0,1] \right) \cap C^{1,2}\left( \left( 0,\|u_0\|_\infty^{-1} \right)\times [0,1] \right)
\end{equation} 
and
\begin{equation}
\label{encad}
\sigma(t) \le u(t,x) \le \Sigma(t)\,, \qquad (t,x) \in \left( 0, \|u_0\|_\infty^{-1} \right)\times (0,1)\,.
\end{equation}
If $v$ denotes the unique solution to $-\partial_x^2 v = u-M$ in $\left( 0,\|u_0\|_\infty^{-1} \right)\times (0,1)$ with homogeneous Neumann boundary conditions and mean value zero, it readily follows from \eqref{aux}, \eqref{auxb}, \eqref{cvg}, and \eqref{encad} that $(u,v)$ is a classical solution to \eqref{he1}-\eqref{he3} with the expected diffusion coefficient $a(u)$. Next, either $\|u(t)\|_\infty \longrightarrow \infty$ as $t\to \|u_0\|_\infty^{-1}$ and the existence proof is complete, or there is a sequence $(t_i)_{i\ge 1}$ such that $t_i\to \|u_0\|_\infty^{-1}$ and $(\|u(t_i)\|_\infty)_{i\ge 1}$ is bounded. In the latter, we carry on the construction of $u$ starting from $u(t_{i_0})$ for a suitable choice of $i_0\ge 1$.

Concerning the uniqueness issue, we shall use the Wasserstein distance. Consider an initial condition $u_0$ fulfilling the assumptions of Proposition~\ref{local2} and let $u_1$ and $u_2$ be two classical solutions to \eqref{he1}-\eqref{he3} enjoying the properties listed in Proposition~\ref{local2} on their respective maximal existence time intervals $\left[ 0, T_{max,1} \right)$ and $\left[ 0, T_{max,2} \right)$. Thanks to the smoothness and positivity of $u_1$ and $u_2$, the transformation defined in Section~\ref{cv} is well-defined for both of them. Consequently, according to Proposition~\ref{altform},  we may associate to $u_i$ a solution $f_i$ to 
\begin{eqnarray*}
\partial_t f_i & = & \partial_y^2(\Psi(f_i))-1+Mf_i\,, \qquad (t,y)\in (0,T_{max,i})\times (0,M)\,,\\
& & \partial_y f_i(t,0)=\partial_y f_i(t,M)=0\,, \qquad t\in (0,T_{max,i})\,,\\
f_i(0,y) & = & f_0(y)\,, \qquad y\in (0,M)\,,
\end{eqnarray*}
for $i=1,2$, the functions $\Psi$ and $f_0$ being defined in \eqref{Psi} and Proposition~\ref{altform}, respectively. Now, if $T\in \left( 0 , \min{\left\{ T_{max,1} , T_{max,2} \right\}} \right)$, $f_1$ and $f_2$ are smooth functions which are bounded from above and from below by positive constants (depending on $T$) and standard arguments ensure that 
$$
\|f_1(t)-f_2(t)\|_1 \le \|f_1(0)-f_2(0)\|_1\ e^{Mt} = 0 \;\;\mbox{ for }\;\; t\in [0,T]
$$
(see, e.g., \cite[Theorem~11.2]{Va07}). Consequently, $f_1=f_2$ and also $F_1=F_2$ in $(0,T)\times [0,M]$, where $F_i$ denotes the indefinite integral of $f_i$ vanishing at zero as in Section~\ref{cv}. This readily implies that $u_1=u_2$ in $(0,T)\times [0,1]$ and the proof of Proposition~\ref{local2} is complete. \qed

\medskip

\begin{Rem} \label{rempos}
In the following, we will always assume that there is $m_0\in (0,M)$ such that $u_0\ge m_0$, though it is not needed to apply Proposition~\ref{local1} when $a\in C^1([0,\infty))$. This is not a restriction since, if $a\in C^1([0,\infty))$, Proposition~\ref{local1} warrants that $u(t)$ is positive and continuous for each $t\in (0,T_{max})$. Instead of $u_0$, we take as an initial condition $u(t_0)$ for some $t_0\in (0,T_{max})$ which satisfies the above positivity requirement. The corresponding solution to \eqref{he1}-\eqref{he3} is then $(t,y)\mapsto u(t_0+t,y)$ for $(t,y)\in (0,T_{max}-t_0)\times [0,1]$ and $T_{max}$ and $T_{max}-t_0$ are simultaneously finite or infinite.
\end{Rem}

%%%%%%%%%%%%%%%%%%%%%%%%%%%%%%%%%%%%%%%%%%%%%%%
%%%%%%%%%%%%%%%%%%%%%%%%%%%%%%%%%%%%%%%%%%%%%%%
\section{Global existence.}\label{glex}
%%%%%%%%%%%%%%%%%%%%%%%%%%%%%%%%%%%%%%%%%%%%%%%
%%%%%%%%%%%%%%%%%%%%%%%%%%%%%%%%%%%%%%%%%%%%%%%

As mentioned in the introduction we first establish that the diffusion coefficient $a(r)=1/r$ is not critical in the sense that all solutions to \eqref{he1}-\eqref{he3} are global in that case (Theorem~\ref{glo}). We in fact prove that the corresponding solution $f$ to the transformed equation \eqref{main} cannot reach zero in finite time: to this end, we deduce an $L^\infty$-bound on $\log{f}$ from a natural Liapunov functional associated to \eqref{main}. In a second step, we extend our analysis to a wider class of diffusion coefficients $a$ (Theorem~\ref{glo1}).

\begin{theo}\label{glo}
Assume that $a(r)=1/r$ and consider a non-negative function $u_0\in C([0,1])$ such that $u_0\ge m_0>0$ and $\langle u_0\rangle=M$ for some $M>0$ and $m_0\in (0,M)$. If $(u,v)$ denotes the corresponding solution to \eqref{he1}-\eqref{he3}, then $T_{max}=\infty$, the trajectory $\{ (u(t),v(t)) \ : \ t\ge 0\}$ is bounded in $L^\infty(0,1;\RR^2)$, and there is $\mu>0$ such that $u(t,x)\ge \mu$ for all $(t,x)\in [0,\infty)\times [0,1]$.
\end{theo}

Global existence is actually true for a wider class of nonlinear diffusion coefficients $a$ such as 
$$
a(r)=\frac{1}{(1+r)^p}\,, \quad p\in (-\infty,1]\,, \;\;\mbox{ and }\;\; a(r) = \frac{1}{(1+r) (\log{(1+r)})^\alpha}\,, \quad \alpha\in (-\infty,1]\,.
$$
Indeed, the above examples fulfil the assumptions of the following result:

\begin{theo}\label{glo1}
Assume that $a\in C^1((0,\infty))$ is a positive function such that $a\not\in L^1(1,\infty)$ and, for each $\varepsilon\in (0,\infty)$, there is $\kappa_\varepsilon>0$ for which
\begin{equation}
\label{gex1}
a(r) \le \varepsilon\ r a(r) + \frac{\kappa_\varepsilon}{r} \;\;\mbox{ for }\;\; r\in (0,1)\,.
\end{equation}
Consider an initial condition $u_0\in C([0,1])$ such that $u_0\ge m_0>0$ and $\langle u_0\rangle=M$ for some $M>0$ and $m_0\in (0,M)$. If $(u,v)$ denotes the corresponding solution to \eqref{he1}-\eqref{he3} given by Proposition~\ref{local2}, then $T_{max}=\infty$ and the trajectory $\{ (u(t),v(t)) \ : \ t\ge 0\}$ is bounded in $L^\infty(0,1;\RR^2)$. In addition, if $a\not\in L^1(0,1)$, there is $\mu>0$ such that $u(t,x)\ge \mu$ for all $(t,x)\in [0,\infty)\times [0,1]$.
\end{theo}

Though the proofs of Theorems~\ref{glo} and~\ref{glo1} follows the same steps, the generality of the latter brings additional technicalities in the proof. We will thus first present the proof of Theorem~\ref{glo} to illustrate the main ideas and then proceed with the proof of Theorem~\ref{glo1}.

\medskip

\noindent{\bf Proof of Theorem~\ref{glo}.} According to Proposition~\ref{local2}, it suffices to provide a control on the $L^\infty$-norm of $u$ on any finite time interval to ensure global existence. By Proposition~\ref{altform} and Remark~\ref{td} such a control can be obtained by showing that the corresponding solution $f$ to \eqref{main}-\eqref{init} is bounded away from zero on any finite time interval. To achieve the latter, we need some more steps. Recalling that $\Psi(r)=\log r$ in that case, we start with the following observation.
\begin{Le}\label{Liap}
The function 
\begin{equation}\label{Lia}
L(t):= \int_0^M\left(\frac{1}{2}|\partial_y\log f(t,y)|^2+\log f(t,y)\right)\ dy\,, \qquad t\in [0,T_{max})\,,
\end{equation}
is a non-increasing function of time.
\end{Le}    
\noindent{\bf Proof of Lemma~\ref{Liap}.} Multiplying \eqref{main} by $\partial_t \log{f}$ and integrating with respect to space over $(0,M)$ we arrive at
\[
\int_0^M\frac{|\partial_t f(t,y)|^2}{f(t,y)}\ dy=-\int_0^M\partial_y\log f(t,y) \partial_t\partial_y\log f(t,y)\  dy - \int_0^M \partial_t\log f(t,y)\ dy + M\ \int_0^M \partial_t f(t,y)\ dy.
\]
Owing to the conservation of mass \eqref{mas}, the last term of the right-hand side of the above identity vanishes and we end up with
\begin{equation}\label{Li}
\frac{dL}{dt}(t)+\int_0^M f(t,y) |\partial_t\log f(t,y)|^2\ dy=0,
\end{equation}
from which Lemma~\ref{Liap} readily follows. \qed

\medskip

At this point, we notice that, if a bound from below turns out to be available for $L$, we might expect to deduce from it an estimate on $\log{f}$ in $L^\infty(0,T;H^1(0,M))$. Thanks to the continuous embedding of $H^1(0,M)$ in $L^\infty(0,M)$, such an estimate will in turn provide an $L^\infty$-estimate on $\log{f}$ from which clearly follows that $f$ is bounded away from zero. The next step in the proof of Theorem~\ref{glo} thus requires a detailed study of $L$. To this end, we define the functional 
\begin{equation}\label{E}
E(g):=\int_0^M\left(\frac{1}{2}|\partial_y g(y)|^2+g(y)\right)\ dy \;\;\mbox{ for }\;\; g\in H^1(0,M)\,,
\end{equation}
so that Lemma~\ref{Liap} states that $t\longmapsto E(\log{f(t)})$ is a non-increasing function of time. We next proceed to find a lower bound on $E$.
\begin{pr}\label{2.3}
If $g\in H^1(0,M)$ satisfies $\langle e^g\rangle=1/M$ then
\begin{equation}\label{2.12}
E(g)\geq \frac{1}{4}\|\partial_y g\|_2^2 -M\log M -M^3\,,
\end{equation}
and
\begin{equation}\label{2.13}
\|g\|_1\leq 2 +M\log M + M^{3/2} \|\partial_y g\|_2\,.
\end{equation}
\end{pr}

\noindent{\bf Proof of Proposition~\ref{2.3}.} We consider $g\in H^1(0,M)$ such that $\langle e^g\rangle=1/M$. Thanks to the continuous embedding of $H^1(0,M)$ in $C([0,M])$, $e^g$ belongs to $C([0,M])$ and the constraint $\langle e^g\rangle=1/M$ and the mean value theorem ensure that there exists $Y_g\in [0,M]$ such that $e^{g(Y_g)}=1/M$. Then 
\begin{eqnarray}
\int_0^M g(y) dy & = & \int_0^M (g(y)-g(Y_g))\ dy + M g(Y_g) \nonumber\\
& \ge & \int_0^M \int_{Y_g}^y \partial_y g(z)\ dzdy - M \ \log{M} \nonumber\\
& \ge & - M^{3/2} \|\partial_y g\|_2 - M \log{M}\,. \label{2.125}
\end{eqnarray}
Thanks to \eqref{2.125} we can estimate $E(g)$ from below and obtain
\begin{eqnarray*}
E(g) & \ge & \frac{1}{4}\ \|\partial_y g\|_2^2 + \frac{1}{4}\ \|\partial_y g\|_2^2 - M^{3/2} \|\partial_y g\|_2 + M^3 - M^3 - M  \log{M} \\
& \ge & \frac{1}{4}\ \|\partial_y g\|_2^2 + \frac{1}{4}\ \left( \|\partial_y g\|_2 - 2 M^{3/2} \right)^2 - M^3 - M \log{M}\,,
\end{eqnarray*}
hence \eqref{2.12}. It also follows from \eqref{2.125} that 
\begin{eqnarray*}
\|g\|_1 & = & 2\int_0^M \max\{0,g(y)\}\ dy- \int_0^M g(y)\ dy \\
& \leq & 2\int_0^M \left( e^{\max\{0,g(y)\}} - 1 \right)\ dy +M\log M + M^{3/2}\|\partial_y g\|_2 \\
& \leq & 2 \int_0^M e^g(y)\ dy +M\log M + M^{3/2} \|\partial_y g\|_2\,,
\end{eqnarray*}
from which \eqref{2.13} follows as $\langle e^g\rangle=1/M$. \qed

We now complete the proof of Theorem~\ref{glo} by noticing that Proposition~\ref{2.3} implies that $\log f\in L^\infty(0,T;H^1(0,M))$ for any $T\in (0,T_{max})$ with a bound that depends neither on $T$ nor on $T_{max}$.  Indeed, since $\langle e^{\log{f}}\rangle= 1/M$ by \eqref{mas}, we infer from \eqref{Lia} and \eqref{2.12} that, for $t\in [0,T_{max})$, 
\[
\frac{1}{4} \|\partial_y\log f(t)\|_2^2\leq E(\log f(t)) +M\log M +M^3 \leq E(\log f_0) +M\log M +M^3.
\]
In view of \eqref{2.13}, the previous bound entails that $\log{f}$ belongs to $L^\infty(0,T_{max};L^1(0,M))$ and we finally conclude with the help of the Poincar\'e-Wirtinger inequality that
\begin{equation}\label{ost}
\|\log f(t)\|_{H^1(0,M)}\leq C(M)(1+E(\log f_0)) \;\;\mbox{ for }\;\; t\in [0,T_{max}).
\end{equation}
Owing to the continuous embedding of $H^1(0,M)$ in $L^\infty(0,M)$ and \eqref{ost}, we conclude that there are
constants $0<\alpha<\beta$ such that $\alpha<f(t,y)<\beta$ for $(t,y)\in [0,T_{max})\times [0,M]$. This property clearly warrants that $T_{max}=\infty$ and ends the proof of Theorem~\ref{glo} for $a(r)=1/r$. \qed

\medskip

\noindent\textbf{Proof of Theorem~\ref{glo1}.} We follow the same strategy as in the previous case and aim at showing that a similar argument prevents the solution $f$ to \eqref{main}-\eqref{init} to hit zero in finite time, the function $\Psi$ being given by \eqref{Psi} (see Section~\ref{cv}). Owing to the positivity and the non-integrability over $(1,\infty)$ of $a$, we note that $\Psi$ is an increasing function from $(0,\infty)$ onto its range with inverse $\Psi^{-1}$ and
\begin{equation}
\label{gex1b}
\Psi \;\;\mbox{ maps}\;\; (0,1) \;\;\mbox{ onto }\;\; (-\infty,0) \;\;\mbox{ with }\;\; \lim_{r\to 0} \Psi(r)=-\infty\,.
\end{equation}
We next report the analogue of Lemma~\ref{Liap}.
\begin{Le}\label{Liapprim}
The function 
$$
L_1(t):= \frac{1}{2} \int_0^M \left| \partial_y\Psi(f(t,y)) \right|^2\ dy+\int_0^M \left( \Psi(f(t,y)) -M\ \Psi_1(f(t,y)) \right)\ dy
$$
is a non-increasing function of time on $[0,T_{max})$, the function $\Psi_1$ being defined by 
\begin{equation}
\label{gex3}
\Psi_1(1):=0 \;\;\mbox{ and }\;\; \Psi_1'(r) := r \Psi'(r) = \frac{1}{r}\ a \left( \frac{1}{r} \right) \,, \qquad r\in (0,\infty)\,.
\end{equation}
\end{Le} 
\noindent\textbf{Proof of Lemma~\ref{Liapprim}.} We multiply \eqref{main} by $\partial_t \Psi(f)$ and integrate with respect to space over $(0,M)$ to obtain
\begin{eqnarray*}
\int_0^M \Psi'(f)\ |\partial_t f|^2\ dy & = & -\int_0^M \partial_y \Psi(f) \partial_t \partial_y \Psi(f)\ dy\\
& & -\int_0^M \partial_t \Psi(f)\ dy +M \int_0^M f\ \Psi'(f)\ \partial_t f\ dy\,,
\end{eqnarray*}
hence
\[
\frac{dL_1}{dt}(t)+\int_0^M \Psi'(f)\ |\partial_t f|^2\ dy=0\,,
\]
and Lemma~\ref{Liapprim} follows in view of the positivity of $\Psi'$. \qed

\medskip

Observe next that, owing to \eqref{gex1} with $\varepsilon=1/M$, we have 
$$
M\ \Psi_1'(s) = \frac{M}{s}\ a\left( \frac{1}{s} \right) \le \frac{M}{s}\ \left( \frac{1}{Ms}\ a\left( \frac{1}{s} \right) + \kappa_{1/M}\ s \right) \le \Psi'(s) + M\ \kappa_{1/M}, \quad s\in (1,\infty)\,,
$$
so that, after integration over $(1,r)$, $r>1$, $M \Psi_1(r) \le \Psi(r) + M \kappa_{1/M}\ r$. Consequently, since $\Psi_1\le 0$ on $(0,1]$, we have
\begin{eqnarray*}
-M\ \int_0^M \Psi_1(f)(t,y)\ dy & \ge & -M\ \int_0^M \mathbf{1}_{[1,\infty)}(f(t,y))\ \Psi_1(f(t,y))\ dy \\
& \ge & - \int_0^M \mathbf{1}_{[1,\infty)}(f(t,y))\ \left( \Psi(f(t,y)) + M \kappa_{1/M}\ f(t,y) \right)\ dy \\
& \ge & -\int_0^M \mathbf{1}_{[1,\infty)}(f(t,y))\ \Psi(f(t,y))\ dy - M \kappa_{1/M}\,,
\end{eqnarray*}
where we have used \eqref{mas} to obtain the last inequality. Now, recalling \eqref{gex1b}, we find
\begin{eqnarray}
L_1(t) & \ge & \frac{1}{2} \int_0^M \left| \partial_y\Psi(f(t,y)) \right|^2\ dy+\int_0^M \mathbf{1}_{(0,1)}(f(t,y))\ \Psi(f(t,y))\ dy - M \kappa_{1/M} \nonumber\\
& \ge & \frac{1}{2} \int_0^M \left| \partial_y\Psi(f(t,y)) \right|^2\ dy+\int_0^M \mathbf{1}_{(-\infty,0)}(\Psi(f(t,y)))\ \Psi(f(t,y))\ dy - M \kappa_{1/M} \nonumber\\
& \ge & E_1\left( \Psi(f(t)) \right) - M\ \kappa_{1/M}\,, \label{gex4}
\end{eqnarray}
with 
$$
E_1(h) := \frac{1}{2} \|\partial_y h\|_2^2 + \int_0^M \mathbf{1}_{(-\infty,0)}(h(y))\ h(y)\ dy\,, \quad h\in H^1(0,M)\,.
$$
The next step is then to find a lower bound on $E_1$.

\begin{Le}\label{lab1}
For $M>0$, we have
\begin{equation}\label{gex5}
E_1(h) \ge \frac{1}{4}\ \|\partial_y h\|_2^2 - M^3 - M \left| \Psi\left( \frac{1}{M} \right) \right|\,,
\end{equation}
and
\begin{equation}\label{gex6}
\|h\|_1\leq M^{3/2} \|\partial_y h\|_2 + M \left| \Psi\left( \frac{1}{M} \right) \right|
\end{equation} 
for every $h\in H^1(0,M)$ satisfying
\begin{equation}
\label{gex7}
\int_0^M \Psi^{-1}(h)(y)\ dy = 1\,.
\end{equation}
\end{Le}
\noindent{\bf Proof of Lemma~\ref{lab1}.} Consider $h\in H^1(0,M)$ enjoying the property \eqref{gex7}. Arguing as in the proof of Proposition~\ref{2.3} , we infer from \eqref{gex7} and the mean value theorem that there is at least one point $Y_h\in [0,M]$ (depending on $h$) such that $\Psi^{-1}(h(Y_h))=1/M$. Consequently,
\begin{eqnarray}
\int_0^M \min{\{0,h(y)\}}\ dy & = & \int_0^M (\min{\{0,h(y)\}}-\min{\{0,h(Y_h)\}})\ dy + M \min{\{0,h(Y_h)\}} \nonumber\\
& \ge & \int_0^M \int_{Y_h}^y \mathbf{1}_{(-\infty,0)}(h(z))\ \partial_y h(z)\ dzdy - M \left| \Psi\left( \frac{1}{M} \right) \right| \nonumber\\
& \ge & - M^{3/2} \|\partial_y h\|_2 - M \left| \Psi\left( \frac{1}{M} \right) \right|\,. \label{2.1275}
\end{eqnarray}
Thanks to \eqref{2.1275} we can estimate $E_1(h)$ from below and obtain
$$
E_1(h) \ge \frac{1}{4}\ \|\partial_y h\|_2^2 + \frac{1}{4}\ \left( \|\partial_y h\|_2 - 2 M^{3/2} \right)^2 - M^3 - M \left| \Psi\left( \frac{1}{M} \right) \right|\,,
$$
hence \eqref{gex5}. We next proceed as in the proof of \eqref{2.1275} and compute
\begin{eqnarray*}
\int_0^M |h(y)| \ dy & = & \int_0^M (|h(y)|-|h(Y_h)|)\ dy + M |h(Y_h)| \\
& \le & \int_0^M \int_{Y_h}^y \mbox{sign}(h(z))\ \partial_y h(z)\ dzdy + M \left| \Psi\left( \frac{1}{M} \right) \right| \\
& \le & M^{3/2} \|\partial_y h\|_2 + M \left| \Psi\left( \frac{1}{M} \right) \right|\,,
\end{eqnarray*}
thus establishing \eqref{gex6}. \qed

\medskip

Owing to \eqref{mas}, $\Psi(f)$ fulfils the property \eqref{gex7} and we infer from Lemma~\ref{Liapprim}, \eqref{gex4}, and \eqref{gex5} that, for $t\in [0,T_{max})$,
\begin{eqnarray*}
\|\partial_y \Psi(f(t))\|_2^2 & \le & 4 E_1(\Psi(f(t))) + 4M^3 + 4M \left| \Psi\left( \frac{1}{M} \right) \right| \\
& \le & 4 L_1(t) + 4M \kappa_{1/M} + C(a,M) \\
& \le & 4 L_1(0) + C(a,M)\,.
\end{eqnarray*}
Moreover, combining the previous inequality with \eqref{gex6}, we obtain that
$$
\|\Psi(f(t))\|_1 \le M^{3/2}\ \left( 4 L_1(0) + C(a,M) \right)^{1/2} + M\ \left| \Psi\left( \frac{1}{M} \right) \right|\,, \quad t\in [0,T_{max})\,.
$$
Thanks to the Poincar\'e-Wirtinger inequality the above two estimates imply that $\{ \Psi(f)(t)\ : \ t\in [0,T_{max})\}$ is bounded in $L^\infty(0,T_{max};H^1(0,M))$. Owing to the continuous embedding of $H^1(0,M)$ in $C([0,M])$, we conclude that there is $K>0$ such that $-K \le \Psi(f)(t,y) \le K$ for $(t,y)\in [0,T_{max})\times [0,M]$. Since $\Psi$ is an increasing function, we actually have $\Psi^{-1}(-K) \le f(t,y) \le \Psi^{-1}(K)$ for $(t,y)\in [0,T_{max})\times [0,M]$ and $\Psi^{-1}(-K)>0$ by \eqref{gex1}. Recalling Remark~\ref{td}, we thus have shown that $T_{max}=\infty$. Finally, since the lower bound on $f$ do not depend on time, the boundedness of the trajectory  $\{ (u(t),v(t)) \ : \ t\ge 0\}$ in $L^\infty(0,1;\RR^2)$ readily follows. Under the additional assumption $a\not\in L^1(0,1)$, the function $\Psi$ maps $(1,\infty)$ onto $(0,\infty)$ and the time-independent upper bound on $f$ provides a positive lower bound on $u$ thanks to \eqref{odw}. \qed

%%%%%%%%%%%%%%%%%%%%%%%%%%%%%%%%%%%%%%%%%%%%%%%
%%%%%%%%%%%%%%%%%%%%%%%%%%%%%%%%%%%%%%%%%%%%%%%
\section{Finite time blow-up.}\label{ftbu}
%%%%%%%%%%%%%%%%%%%%%%%%%%%%%%%%%%%%%%%%%%%%%%%
%%%%%%%%%%%%%%%%%%%%%%%%%%%%%%%%%%%%%%%%%%%%%%%

In this section, we return to the question of finite time blow-up and show that it can be proved for a large class of diffusion coefficients $a$ with the help of a suitable virial identity. Besides providing finite time blow-up results for functions $a$ which were not included in \cite{CW08}, it also provides an alternative and simpler proof of \cite[Theorem 3.2]{CW08}. 

\begin{theo}\label{blowup}
Let $a\in C^1((0,\infty))$ be a positive function such that $a\in L^1(1,\infty)$ and there is a concave function $B$ for which
\begin{eqnarray}
0 \le -r A(r) & \le & B(r) \;\;\mbox{ with }\;\; A(r) = - \int_r^\infty a(s)\ ds\,, \qquad r\in (0,\infty)\,, \label{bu1} \\
\lim_{r\to \infty} \frac{B(r)}{r} & = & 0\,. \label{bu2}
\end{eqnarray}
Then, for any $M>0$, there is at least an initial condition $u_0\in C([0,1])$ satisfying $u_0\ge m_0>0$ and $\langle u_0\rangle=M$ for some $m_0\in (0,M)$ for which the first component $u$ of the corresponding solution to \eqref{he1}-\eqref{he3} blows up in finite time, i.e. $T_{max}<\infty$.
\end{theo}

Observe that Theorem~\ref{blowup} applies when $a(r)\le C (1+r)^{-p}$ for some $C>0$ and $p\in (1,2]$ (with $B(r)=C r (1+r)^{1-p}/(p-1)$). In that case, we thus recover the blow-up results from \cite{CW08} by a different method. However, Theorem~\ref{blowup} also applies if there exist $C>0$ and $\alpha\in (1,2]$ such that $a(r)\leq C(1+r)^{-1}\left(\log(1+r)\right)^{-\alpha}$ (with $B(r)=C r (\log{(1+r)})^{1-\alpha} / (\alpha-1) $). 

\noindent{\bf Proof of Theorem \ref{blowup}.} Let $q>2$. Introducing $L_q(t):=\|U(t)\|_q^q/q$ (recall that $U(t)$ is the indefinite integral of $u(t)$ vanishing at zero, see Section~\ref{cv}), we infer from \eqref{B}, \eqref{bu1}, the fact that $U(t,0)=M-U(t,1)=0$, and the negativity of $A$ that 
\begin{eqnarray}
\frac{dL_q}{dt} & = & \int_0^1 U^{q-1}\ \partial_xA(u)\  dx +\int_0^1 u\ (U-Mx)\ U^{q-1}\ dx \nonumber\\
& = & \left[ U^{q-1}\ A(u) \right]_{x=0}^{x=1} - (q-1)\ \int_0^1 U^{q-2}\ u\ A(u)\ dx + \left[ \frac{U^{q+1}}{q+1} \right]_{x=0}^{x=1} - M\ \left[ x\ \frac{U^q}{q} \right]_{x=0}^{x=1} + M\ L_q \nonumber\\
& \le & (q-1)\ \int_0^1 U^{q-2}\ B(u)\ dx + M\ L_q - \frac{M^{q+1}}{q(q+1)}\,. \label{bu3}
\end{eqnarray}
By Jensen's inequality with measure $B(u)\ dx$, we have
\begin{eqnarray*}
\int_0^1 U^{q-2}\ B(u)\ dx & = & \left( \int_0^1 B(u)\ dz \right)\ \int_0^1 \left( U^q \right)^{(q-2)/q}\ \frac{B(u)\ dx}{\int_0^1 B(u)\ dz}\\
& \le & \left( \int_0^1 B(u)\ dx \right)^{2/q}\ \left( \int_0^1 U^q\ B(u)\ dx \right)^{(q-2)/q}\\
& \le & \left( q L_q \right)^{(q-2)/q}\ \left( \int_0^1 B(u)\ dx \right)^{2/q}\ \left( \int_0^1 B(u)\ U^q\ \frac{dx}{qL_q} \right)^{(q-2)/q}\,.
\end{eqnarray*}
We now use the concavity of $B$, the Jensen inequality (with measure $dx$ and $U^q\ dx$), and the conservation of mass \eqref{swing} to obtain
\begin{eqnarray*}
\int_0^1 U^{q-2}\ B(u)\ dx & \le & \left( q L_q \right)^{(q-2)/q}\ \left[ B\left( \int_0^1 u\ dx \right) \right]^{2/q}\ \left[ B\left( \int_0^1 u\ U^q\ \frac{dx}{qL_q} \right) \right]^{(q-2)/q}\\
& \le & \left( q L_q \right)^{(q-2)/q}\ \left[ B(M) \right]^{2/q}\ \left[ B\left( \frac{M^{q+1}}{q(q+1)L_q} \right) \right]^{(q-2)/q}\\
& \le & \left[ B(M) \right]^{2/q}\ \left( \frac{M^{q+1}}{q+1} \right)^{(q-2)/q}\ \left[ \beta\left( \frac{M^{q+1}}{q(q+1)L_q} \right) \right]^{(q-2)/q}\,,
\end{eqnarray*}
where $\beta(r)=B(r)/r$ for $r>0$. Recalling \eqref{bu3} and introducing 
$$
\Lambda_M(r) := M\ r + (q-1)\ \left[ B(M) \right]^{2/q}\ \left( \frac{M^{q+1}}{q+1} \right)^{(q-2)/q}\ \left[ \beta\left( \frac{M^{q+1}}{q(q+1) r} \right) \right]^{(q-2)/q} - \frac{M^{q+1}}{q(q+1)}
$$
for $r>0$, we conclude that
\begin{equation}
\label{bu4}
\frac{dL_q}{dt}(t) \le \Lambda_M(L_q(t))\,, \qquad t\in [0,T_{max})\,.
\end{equation}
Assume now for contradiction that $T_{max}=\infty$. Since $\Lambda_M(0)<0$ by \eqref{bu2}, there is $\vartheta_M>0$ such that $\Lambda_M(r)<0$ for $r\in [0,\vartheta_M]$. Then, if $L_q(0)<\vartheta_M$, we readily infer from \eqref{bu4} and the negativity of $\Lambda_M$ in $[0,\vartheta_M]$ that $L_q(t)<\vartheta_M$ and $dL_q(t)/dt \le \sup_{[0,\vartheta_M]}{\{\Lambda_M\}}< 0$ for all $t\ge 0$. This property implies that $L_q$ becomes negative at a finite time and thus contradicts the non-negativity of $L_q$. Therefore, $T_{max}<\infty$ provided that $L_q(0)<\vartheta_M$. The latter condition is fulfilled as soon as $u_0$ is sufficiently concentrated near $x=1$, whatever the value of the mass $M$ is, and the proof of Theorem~\ref{blowup} is complete. \qed

%%%%%%%%%%%%%%%%%%%%%%%%%%%%%%%%%%%%%%%%%%%%%%%
%
\textbf{Acknowledgement.} This paper was prepared during T.~Cie\'slak's one-month visit at the Institut de Math\'ematiques de Toulouse, Universit\'e Paul Sabatier. T.~Cie\'slak would like to express his gratitude for the invitation, support, and hospitality. 
%
%%%%%%%%%%%%%%%%%%%%%%%%%%%%%%%%%%%%%%%%%%%%%%%

%%%%%%%%%%%%%%%%%%%%%%%%%%%%%%%%%%%%%%%%%%%%%%%
%%%%%%%%%%%%%%%%%%%%%%%%%%%%%%%%%%%%%%%%%%%%%%%

%%%%%%%%%%%%%%%%%%%%%%%%%%%%%%%%%%%%%%%%%%%%%%%
%%%%%%%%%%%%%%%%%%%%%%%%%%%%%%%%%%%%%%%%%%%%%%%


\begin{thebibliography}{20}

\bibitem{BN94}
P.~Biler and T.~Nadzieja, \textit{Existence and nonexistence of solutions for a model of gravitational interaction of particles~I}, Colloq. Math. {\bf 66} (1994), 319--334.

\bibitem{BCLxx}
A.~Blanchet, J.A.~Carrillo, and Ph.~Lauren\c{c}ot,
\textit{Critical mass for a Patlak-Keller-Segel model with degenerate diffusion in higher dimensions}, Calc. Var. Partial Differential Equations, to appear.

\bibitem{CF05}
J.A.~Carrillo and K.~Fellner, \textit{Long-time asymptotics via entropy methods for diffusion dominated equations}, Asymptot. Anal. {\bf 42} (2005), 29--54.

\bibitem{CGT04}
J.A.~Carrillo, M.P.~Gualdani, and G.~Toscani, \textit{Finite speed propagation in porous media by mass transportation methods}, C. R. Math. Acad. Sci. Paris \textbf{338} (2004), 815--818.

\bibitem{Ch03}
P.-H.~Chavanis, \textit{Generalized thermodynamics and Fokker-Planck equations. Applications to stellar dynamics and two-dimensional turbulence}, Phys. Rev. E \textbf{68} (2003), 036108.

\bibitem{ChS04} P.-H.~Chavanis and C.~Sire, \textit{Anomalous diffusion and collapse of self-gravitating Langevin particles in $D$ dimensions}, Phys. Rev. E \textbf{69} (2004), 016116.

\bibitem{cl1}
T.~Cie\'slak and Ph.~Lauren\c{c}ot, \textit{Finite time blow-up for radially symmetric solutions to a critical quasilinear Smoluchowski-Poisson system}, C. R. Acad. Sci. Paris, S\'er.~I, to appear.

\bibitem{cl2}
T.~Cie\'slak and Ph.~Lauren\c{c}ot, \textit{Finite time blow-up for a one-dimensional quasilinear parabolic-parabolic chemotaxis system}, (submitted).

\bibitem{CW08}
T.~Cie\'slak and M.~Winkler, \textit{Finite time blow-up in a quasilinear system of chemotaxis}, Nonlinearity \textbf{21} (2008), 1057-1076.

\bibitem{DFW08}
M.~Di Francesco and M.~Wunsch, \textit{Large time behavior in Wasserstein spaces and relative entropy for bipolar 
drift-diffusion-Poisson models}, Monatsh. Math. {\bf 154} (2008), 39--50.

\bibitem{DP04}
J.~Dolbeault and B.~Perthame, \textit{Optimal critical mass in the two-dimensional Keller-Segel model in $\RR^2$}, C. R. Math. Acad. Sci. Paris \textbf{339} (2004), 611--616.

\bibitem{HV97}
M.A.~Herrero and J.J.L.~Vel\'azquez, \textit{A blow-up mechanism for a chemotaxis model},  Ann. Scuola Norm. Super. Pisa Cl. Sci. \textbf{24} (1997), 633--683.

\bibitem{GZ98}
H.~Gajewski and K.~Zacharias, \textit{Global behavior of a reaction-diffusion system modelling chemotaxis},
Math. Nachr. {\bf 195} (1998), 77--114. 

\bibitem{JL92}
W.~J\"{a}ger and S.~Luckhaus, \textit{On explosions of solutions to a system of partial differential equations modelling chemotaxis}, Trans. Amer. Math. Soc. \textbf{329} (1992), 819--824.

\bibitem{KS70}
E.F.~Keller and L.A.~Segel, \textit{Initiation of slime mold aggregation viewed as an instability}, J. Theoret. Biol. {\bf 26} (1970), 399--415.

\bibitem{LSU} O.A.~Ladyzhenskaya, V.A.~Solonnikov, and N.~Uraltseva, \textit{Linear and Quasilinear Equations of Parabolic Type}, Amer. Math. Soc., Providence, 1968.

\bibitem{Na95}
T.~Nagai, \textit{Blow-up of radially symmetric solutions to a chemotaxis system}, Adv. Math. Sci. Appl. {\bf 5} (1995), 581--601.

\bibitem{Su07}
Y.~Sugiyama, \textit{Application of the best constant of the Sobolev inequality to degenerate Keller-Segel models}, Adv. Differential Equations \textbf{12} (2007), 121--144.

\bibitem{Va07}
J.L.~V\'azquez, \textit{The Porous Medium Equation. Mathematical Theory}, Oxford Math. Monographs, Oxford Univ. Press, Oxford, 2007.

\bibitem{Vi03}
C.~Villani, \textit{Topics in Optimal Transportation}, Graduate Studies in Mathematics \textbf{58}, Amer. Math. Soc., Providence, 2003. 

\end{thebibliography}
\end{document}